\documentclass{amsart}
\usepackage{enumerate,amsmath,amssymb,amsfonts,bm}
\usepackage{frcursive,fullpage,graphicx,psfrag}

\newcommand{\RR}{\mathbb{R}} 
\newcommand{\ZZ}{\mathbb{Z}}
\newcommand{\QQ}{\mathbb{Q}}

\newcommand{\EE}{\mathbb{E}}
\newcommand{\NN}{\mathbb{N}}
\newcommand{\LL}{\mathbb{L}}
\newcommand{\cA}{\mathcal{A}}
\newcommand{\cC}{\mathcal{C}}
\newcommand{\cF}{\mathcal{F}}
\newcommand{\cS}{\mathcal{S}}

\newcommand{\cO}{\mathcal{O}}
\newcommand{\cT}{\mathcal{T}}
\newcommand{\cE}{\mathcal{E}}

\newcommand{\ux}{\mathbf{x}}
\newcommand{\uy}{\mathbf{y}}

\begin{document}

\author{Abdelmalek Abdesselam}
\address{Abdelmalek Abdesselam,
Department of Mathematics,
P. O. Box 400137,
University of Virginia,
Charlottesville, VA 22904-4137, USA}
\email{malek@virginia.edu}

\title{QFT, RG, and all that, for mathematicians, in eleven pages}

\begin{abstract}
We present a quick introduction to quantum field theory and Wilson's theory of the renormalization
group from the point of view of mathematical analysis. The presentation is geared primarily towards a
probability theory, harmonic analysis and dynamical systems theory audience.
\end{abstract}

\maketitle

\section{Introduction}
To say that quantum field theory (QFT) has exerted a profound influence on recent mathematical developments is a
banal statement. Ideas from QFT have been shown to be relevant for the understanding of 
knot invariants~\cite{WittenKnots},
four-manifolds~\cite{Witten4manifolds} and questions in enumerative algebraic geometry~\cite{CandelasOGP}.
Not only low-dimensional topology but also high-dimensional topology benefited from QFT (see, e.g.,~\cite{Lurie}). 
The connections between two-dimensional conformal field theory (CFT) and the geometric Langlands correspondence are
well-known~\cite{Gaitsgory,FrenkelB,BeilinsonD,Frenkel}.
Moreover, the latter has recently been shown to be related to higher-dimensional QFT~\cite{KapustinW}. 
It is therefore not surprising that, in recent times,
there has been an increased mathematical interest for QFT.
Several books have appeared where mathematicians took on the task of explaining QFT to other mathematicians (see, e.g.,~\cite{Folland,deFariaM,Costello}).
Yet, it would be fair to say that this interest came mostly from
domains of mathematics such as geometry, topology and representation theory,
while the general area known as analysis has been lagging behind. The proportion of young mathematical analysts
who are working on the foundational questions posed by QFT, compared to other areas of analysis such as partial differential equations,
is very small. Hopefully, the introduction presented here will help make the subject more approachable to such analysts.
The purpose of this article, geared primarily towards a probability theory, harmonic analysis and dynamical systems theory audience,
is to give an idea of some of the key mathematical problems posed by QFT and to explain how
Wilson's renormalization group (RG) theory offers a strategy for solving them. 
The main problem relates to the construction of QFT functional integrals. Unfortunately, this is not addressed in the above-mentioned books.
Indeed, these books only consider the construction in the sense of formal power series.
The main problem will be presented in \S\ref{fundamentalsec}.
Before that, \S\ref{scalinglimsec} will provide a motivation for the study of this problem coming from the analysis of scaling
limits for models in statistical mechanics such as the famous two-dimensional Ising model.
Then in \S\ref{strategysec}, we will provide a rough outline of the RG strategy for solving this problem.
Finally, in \S\ref{hierarchicalsec} we will discuss hierarchical models which constitute a useful testing ground for
rigorous RG methods.
Note that because of the self-imposed page limitation and the intent not to obscure the big picture,
many technical details will be omitted from the discussion. We tried to do so while sacrificing as little mathematical precision as possible.
At the end of the article, we will list references where the rigorous mathematical details may be found.

\section{Scaling limits}\label{scalinglimsec}
A major theme in today's probability theory is the study of scaling limits of models from
statistical mechanics in relation to CFT.
A typical example is that of the Ising model on a two-dimensional lattice
and at the critical temperature.
If $(\sigma_{\ux})_{\ux\in\ZZ^2}$ denotes the random configuration of Ising spins,
one can associate to it a random generalized function or distribution
$\phi_r=L^{\frac{15}{8}r}\sum_{\ux\in\ZZ^2} \sigma_{\ux} \delta_{L^r\ux}$.
Here $L$ is some fixed number greater than $1$ which serves as a yardstick for measuring changes of scale.
In the context of the dyadic decompositions frequently used in harmonic analysis,
one picks $L=2$. The notation $\delta_{L^r\ux}$ refers to the delta function located at $L^r\ux\in\RR^2$.
As for $r\in\ZZ$, it plays the role of an ultraviolet (UV) or short distance cut-off since one can think of the Ising model
as now living on the lattice $(L^r\ZZ)^2$ whose mesh $L^r$ is taken to $0$.
The scaling limit is the generalized random field $\phi_{-\infty}$, or simply $\Phi$, obtained by this construction
when $r\rightarrow -\infty$. The uniqueness and conformal invariance of this scaling limit was shown
in the recent work~\cite{CamiaGN} which builds on~\cite{ChelkakHI}.
The main reason to consider such scaling limits which, by definition, live on the continuum
is that they are universal objects with enhanced symmetry. Instead of lattice symmetries one gets full invariance by translation,
rotation and, here also, by scale transformation. Since interactions are local (e.g., nearest neighbor ones for the Ising model),
one expects these symmetries to hold locally, hence conformal invariance.
The latter was introduced in the present context by Polyakov in~\cite{Polyakov} (see also~\cite{FultonRW} for its early history
in physics).
Important quantities of interest are the moments or correlators
$\EE\left[\Phi(x_1)\ldots\Phi(x_n)\right]$
of $\Phi$. These are expected to be distributions with singular support
on the big diagonal where $x_i=x_j$ for some $i\neq j$.
This relates to the fact such a random field is a generalized one whose sample paths are given by distributions
rather than ordinary functions.
For the Ising model, the $2$-point function $\EE[\Phi(x_1)\Phi(x_2)]$ decays like $|x_1-x_2|^{-\frac{1}{4}}$
which reflects a scaling dimension $[\Phi]=\frac{1}{8}$ for the (elementary) field
$\Phi$.
In two dimensions a free field should behave logarithmically and thus with a scaling dimension $0$. The difference
$\frac{1}{8}-0$ is an example of anomalous dimension.
Via the Schwartz nuclear theorem, the above correlators seen as distributions are obtained from
the moments of honest random variables $\Phi(f)=\int \Phi(x)f(x)\ {\rm d}^{2}x$ where the field is smeared with a test function $f$.
For a dimension of space $d=2$, such variables are the limits when $r\rightarrow -\infty$ of unit lattice quantities
\begin{equation}
\phi_r(f)=\sum_{\ux\in\ZZ^2} \sigma_{\ux}\ L^{(d-[\Phi])r} f(L^{r}x)
\label{spindiluted}
\end{equation}
which involve the diluted test function $L^{(d-[\Phi])r} f(L^{r}\cdot)$.

Another feature which emerges when considering the scaling limit in the continuum is that correlators have a precise asymptotic
expansion in the limit where two of the evaluation points coincide and, moreover, the shape of such an expansion is uniform
with respect to the other points. This is the operator product expansion (OPE) which lies at the foundation of 
CFT in physics. Efforts to mathematically capture this structure can be seen in such frameworks as Borcherds' vertex operator 
algebras~\cite{Borcherds,FrenkelB}, Beilinson and Drinfeld's chiral
algebras~\cite{BeilinsonD,Gaitsgory} or Costello and Gwilliam's factorization algebras~\cite{CostelloG}.
Terms in such expansions may be viewed as mixed moments which, in addition to the original random field, involve suitably defined pointwise
squares, third powers etc. of that field. This generalizes the notion of Wick power for a Gaussian field. If the $2$-point correlation for the
squared field decays with a power which is not twice that of the original field, then one says that the squared field displays an anomalous
dimension of its own. For the Ising model, the successes mentioned earlier (e.g.,~\cite{CamiaGN})
were made possible by very special features of two-dimensional lattice
models: exact solutions~\cite{Wu} and suitable notions
of discrete holomorphic 
functions~\cite{Smirnov}, or the SLE~\cite{Schramm}.
For more general models where such tools are not available, Wilson's Nobel Prize winning
theory~\cite{WilsonII,WilsonF}
of the RG is 
about the only game in town.

\section{The fundamental problem}\label{fundamentalsec}
The origins of the RG come from QFT where, from the mathematical point of view, the fundamental problem is to
give a meaning to and study the properties of expressions such as
\[
\EE\left[\cO_{A_1}(x_1)\cdots\cO_{A_n}(x_n)\right]=
\frac{\int_{\cF}\cO_{A_1}(x_1)\cdots\cO_{A_n}(x_n)\ e^{-S(\Phi)}\ {\rm D}\Phi}{\int_{\cF}e^{-S(\Phi)}\ {\rm D}\Phi}\ .
\]
The integrals are over a space $\cF$ of ``functions'' $\Phi:\RR^d\rightarrow\RR$, with ${\rm D}\Phi$ denoting the ``Lebesgue measure''
on this infinite-dimensional space. For most applications it is enough to take the space of tempered distributions $\cF=S'(\RR^d)$.
As for the functional $S$, a typical example is
\[
S(\Phi)=\int_{\RR^d}\left\{
\frac{1}{2}(\nabla\Phi)^2(x)+\mu\ \Phi(x)^2+g\ \Phi(x)^4
\right\}\ {\rm d}^dx
\]
which corresponds to the so-called $\phi_d^4$ model. Initially, the latter was only thought of as a toy model for more physical ones
such as quantum electrodynamics which describes particles seen in nature (photons, electrons and positrons).
However, this is no longer the case since the discovery of the Higgs particle~\cite{ATLAS,CMS}.
Finally, a local observable $\cO_A(x)$ stands for a function of the field $\Phi$ and its derivatives at the point $x\in\RR^d$
such as $\Phi(x), \Phi(x)^2, \Phi(x)^3\partial_i\Phi(x)$, etc. The different species of such observables are labelled by $A\in\cA$.
One thus avoids precise plethystic notations for ``monomials of monomials''.
The OPE (see, e.g.,~\cite{WittenOPE}) is the asymptotic expansion
\[
\EE\left[\cO_{A_1}(x_1)\cO_{A_2}(x_2)\cO_{A_3}(x_3)\cdots\cO_{A_n}(x_n)\right]=
\sum_{j=0}^{\infty} \cC_j(x_1-x_2)
\EE\left[\cO_{B_j}(x_1)\cO_{A_3}(x_3)\cdots\cO_{A_n}(x_n)\right]
\]
when $x_2\rightarrow x_1$.
The nontrivial requirement is that the functions $\cC_j$ and the (composite) fields $\cO_{B_j}$ should remain the same for any $n$ and
for whatever locations $x_3,\ldots,x_n$ of the spectator fields.
See~\cite{KangM} for a recent probabilistic investigation of the OPE
not far removed from the one presented here.

\section{The RG strategy}\label{strategysec}
The following is a distillation by the author of the ideas of Wilson~\cite{WilsonII,WilsonK}
and Wegner~\cite{Wegner} regarding the RG strategy for solving
the above fundamental problem. First, one combines the kinetic part with the nonexistent Lebesgue measure ${\rm D}\Phi$
and turns them into a Gaussian measure ${\rm d}\mu_{C_{-\infty}}$ with covariance
$C_{-\infty}(x,y)$ or $C_{-\infty}(x-y)$ (slight abuse of notation)
given in Fourier space by $\widehat{C}_{-\infty}(\xi)=|\xi|^{-d+2[\phi]}$.
Traditionally, the scaling dimension $[\phi]$ is taken to be the canonical one, i.e., $\frac{d-2}{2}$;
however, it is important to keep the greater
generality of allowing other values corresponding to fractional powers of the Laplacian.
Then, one mollifies this covariance at distance scale $L^r$ by introducing the cut-off covariance $C_r$ such that
$\widehat{C}_{r}(\xi)=\widehat{C}_{-\infty}(\xi)\eta(L^r|\xi|)$ for some function $\eta$ which is $1$ near $0$ and is
$0$ for arguments greater
than $1$. One also needs a volume cut-off in a box of size $L^s$, $s\rightarrow\infty$, but this will be ignored in the present
discussion.
The log-moment generating function for the wanted random field is
\[
\cS^T(f)=\lim_{r\rightarrow -\infty}
\log\left\{\frac{\int_{S'(\RR^d)}\exp\left(-V_r(\phi_r)+\phi_r(f)\right)\ {\rm d}\mu_{C_r}(\phi_r)}
{\int_{S'(\RR^d)}\exp\left(-V_r(\phi_r)\right)\ {\rm d}\mu_{C_r}(\phi_r)}\right\}
\]
where the bare potentials $V_r$ are given by
\[
V_r(\phi_r)=\int_{\RR^d} \left\{
\mu_r :\phi_r^2:_{C_r}(x)+g_r :\phi_r^4:_{C_r}(x)
\right\}\ {\rm d}^dx\ .
\]
The switch to Wick powers $:\phi_r^2:_{C_r}(x)= \phi_r(x)^2-C_r(0)$ and
$:\phi_r^4:_{C_r}(x)=\phi_r(x)^4-6\ C_r(0)\phi_r(x)^2+3\ C_r(0)^2$ is a simple linear triangular change of parameters.
The input for the construction is the choice of bare ansatz, i.e., the sequence $(g_r,\mu_r)_{r\in\ZZ}$ or rather its germ
at $-\infty$. A preparatory step for the use of the RG is a simple scaling transformation which rewrites
the quantities of interest ``on the unit lattice''. Indeed, the field $\phi_r$ sampled according to the Gaussian measure ${\rm d}\mu_{C_r}$
has the same law as $L^{-[\phi]r}\phi_0(L^{-r}\cdot)$ where $\phi_0$ or simply $\phi$ is sampled according to ${\rm d}\mu_{C_0}$ with
UV cut-off at scale $1$.
As a result
\begin{equation}
\frac{\int_{S'(\RR^d)}\exp\left(-V_r(\phi_r)+\phi_r(f)\right)\ {\rm d}\mu_{C_r}(\phi_r)}
{\int_{S'(\RR^d)}\exp\left(-V_r(\phi_r)\right)\ {\rm d}\mu_{C_r}(\phi_r)}=
\frac{\int_{S'(\RR^d)}\exp\left(-V^{(r,r)}[f](\phi)\right)\ {\rm d}\mu_{C_0}(\phi)}
{\int_{S'(\RR^d)}\exp\left(-V^{(r,r)}[0](\phi)\right)\ {\rm d}\mu_{C_0}(\phi)}
\label{startRG}
\end{equation}
where 
\[
V^{(r,r)}[0](\phi)=\int_{\RR^d} \left\{
\mu^{(r,r)} :\phi^2:_{C_0}(x)+g^{(r,r)} :\phi^4:_{C_0}(x)
\right\}\ {\rm d}^dx
\]
with $g^{(r,r)}=L^{(d-4[\phi])r}g_r$, $\mu^{(r,r)}=L^{(d-2[\phi])r}\mu_r$
and where $V^{(r,r)}[f]$ is the same thing but with the addition of an extra term coming from the test function $f$, namely,
\[
\int_{\RR^d} \phi(x)\ L^{(d-[\phi])r}f(L^r x)\ {\rm d}^{d}x\ .
\]
The latter features the diluted test function in the same way as in (\ref{spindiluted}),
if one extends the spin field $\sigma$ so as to be constant
on unit cells $\Delta=\ux+[0,1)^2$, $\ux\in\ZZ^2$, and if one ignores the variation of this diluted test function on such cells.

Now one exploits the covariance decomposition $C_0=\Gamma+C_1$ (which defines the fluctuation covariance $\Gamma$)
in order to realize $\phi$ as the sum of two independent Gaussian fields $\zeta$ and $\psi$ with respective covariances $\Gamma$ and $C_1$.
Since $\psi$ has the same law as $L^{-[\phi]}\phi(L^{-1}\cdot)$, one can trade $\psi$ for the original field variable $\phi$ and
obtain the identity
\[
\int_{S'(\RR^d)}\exp\left(-V(\phi)\right)\ {\rm d}\mu_{C_0}(\phi)
=\int_{S'(\RR^d)}\exp\left(-\widetilde{V}(\phi)\right)\ {\rm d}\mu_{C_0}(\phi)
\]
where
\[
\widetilde{V}(\phi)=-\log\left\{
\int_{S'(\RR^d)}\exp\left(-V(\zeta+L^{-[\phi]}\phi(L^{-1}\cdot))\right)\ {\rm d}\mu_{\Gamma}(\zeta)
\right\}\ .
\]
The RG map is the transformation $V\rightarrow V'$ where $V'=\widetilde{V}+\delta b(V)$ is obtained by adding a field-independent
quantity $\delta b(V)$ which, for the purposes of the present discussion,
can be thought of as minus the value of $\widetilde{V}(\phi)$ at $\phi=0$.
One application of the RG transformation thus gives an identity
\[
\int_{S'(\RR^d)}\exp\left(-V(\phi)\right)\ {\rm d}\mu_{C_0}(\phi)
=e^{\delta b(V)}\times \int_{S'(\RR^d)}\exp\left(-V'(\phi)\right)\ {\rm d}\mu_{C_0}(\phi)\ .
\]
By repeating this operation for the numerator and denominator in (\ref{startRG}) one generates sequences of
potentials $V^{(r,r)}[f]\rightarrow V^{(r,r+1)}[f]\rightarrow\cdots$ and $V^{(r,r)}[0]\rightarrow V^{(r,r+1)}[0]\rightarrow\cdots$
respectively. In other words, for any $q\in\ZZ$, $V^{(r,q)}[f]=RG^{q-r}\left(V^{(r,r)}[f]\right)$ provided $r$ is sufficiently negative.
The denominator corresponds to the special case $f=0$. 
The $r\rightarrow -\infty$ limit gives an exact formula for the joint cumulant generating function
of the desired random field:
\begin{equation}
\cS^T(f)=\sum_{q\in\ZZ} \left\{
\delta b\left(V^{(-\infty,q)}[f]\right)
-\delta b\left(V^{(-\infty,q)}[0]\right)
\right\}
\label{dbseries}
\end{equation}
where, for any fixed (logarithmic) scale $q\in\ZZ$, $V^{(-\infty,q)}[f]=\lim_{r\rightarrow -\infty} V^{(r,q)}[f]$.
The success of this RG strategy rests on doing two things.

\noindent{\bf 1) Stabilizing the bulk:}
Note that the potentials $V^{(r,q)}[0]$ coming from the partition functions in the denominator
must live in a space $\cE_{\rm bulk}$ of potentials which have coefficients or couplings in front of $:\phi^2:(x)$ and
$:\phi^4:(x)$ (as well as many other new terms generated by the RG) which are uniform in space, i.e., do not depend on the location
$x$. Stabilizing the bulk means controlling the limits which give rise to the sequence
$\cT_{\rm ideal}=\left(V^{(-\infty,q)}[0]\right)_{q\in\ZZ}$ in $\cE_{\rm bulk}$.
This (ideal) sequence is typically made of potentials which are not clean ones consisting of a linear combination of $:\phi^4:$ and $:\phi^2:$
only, such as the original $V^{(r,r)}[0]$. The sequence $\cT_{\rm ideal}$
is also a complete orbit or trajectory for the dynamical system $RG$ acting on $\cE_{\rm bulk}$. It starts at
a UV fixed point $V_{\rm UV}$, when $q\rightarrow -\infty$, and ends at an infrared (IR) fixed point $V_{IR}$, when $q\rightarrow\infty$.
Most models studied in constructive QFT are massive ones where $V_{\rm UV}$ is the massless Gaussian fixed point $V_{\rm Gauss}$
(or ``Gaussian free field'') and $V_{\rm IR}$ is the so-called high temperature fixed point $V_{\rm HT}$ that
corresponds to white noise on $\RR^d$.
The crude approximations $e^x\simeq 1+x$ and $\log(1+x)\simeq x$ allow one to compute the linearization of the map $RG$ at
$V_{\rm Gauss}$. The Wick powers $\int_{\RR^d}:\phi^k:_{C_0}(x)\ {\rm d}^{d}x$ are eigenvectors for this linearization
with eigenvalues $L^{d-k[\phi]}$. For $d>4$ and $[\phi]=\frac{d-2}{2}$ the triviality of the $\phi^4$ model
amounts, in this setting, to the absence of a third fixed point and the fact $V_{\rm Gauss}$ has a one-dimensional unstable manifold
where all corresponding random fields are Gaussian (restriction to even potentials with $k\in 2\NN$ is understood and $k=0$
has been eliminated by the introduction of the $\delta b$ shifts).
A more interesting situation is $d=3$ with $[\phi]=\frac{3-\epsilon}{4}$ where $\epsilon>0$ is a small bifurcation parameter.
This is the Brydges-Mitter-Scoppola (BMS) model~\cite{BMS}
where a new nontrivial fixed point $V_{\rm nontriv}$ appears in the $\int :\phi^4:$ direction
at distance $\epsilon$ from $V_{\rm Gauss}$. The fixed point $V_{\rm nontriv}$ has a codimention one stable manifold, with unstable direction
in the general direction of the mass term $\int :\phi^2:$ leading to $V_{\rm HT}$.
A particular example of bare ansatz consists in choosing $g^{(r,r)}=g$ and $\mu^{(r,r)}=\mu$ fixed. Similarly to tuning the temperature
in the Ising model to its critical value, one can pick a critical value $\mu=\mu_{c}(g)$ so that the ($r$-independent) $V^{(r,r)}[0]$
lie on the stable manifold of $V_{\rm nontriv}$.
This gives rise to a scaling limit similar to that for the critical Ising model discussed at the beginning.
In fact, it is conjectured that a suitable random field obtained along these lines for $d=2$ and $[\phi]=0$ is the same as the random field 
$\Phi$ constructed in~\cite{CamiaGN}. In this case $\cT_{\rm ideal}$ is the constant sequence equal to $V_{\rm nontriv}$.
For the BMS model and this type of bare ansatz, step {\bf 1)} has been done rigorously in~\cite{BMS}.
In~\cite{AbdesselamCMP}, the author constructed a connecting orbit which joins $V_{\rm Gauss}$ to $V_{\rm nontriv}$.
It is not difficult to elaborate on the proof therein in order to produce a bare ansatz which results in a trajectory
$\cT_{\rm ideal}$ which is that connecting
orbit. Thus step {\bf 1)} is essentially solved also for a random field which is not self-similar but has a short distance scaling limit
which is Gaussian and a large distance one which is not and corresponds to the nontrivial fixed point. These two notions of scaling limit,
for a generalized random field in the continuum, are the ones defined, e.g., in~\cite{Dobrushin2}.

\noindent{\bf 2) Controlling the deviations:}
The previous RG map is defined over the bigger (extended) space $\cE_{\rm ext}$ that allows couplings in front of $:\phi^k:(x)$
to depend on the location $x$. The previous space $\cE_{\rm bulk}\subset\cE_{\rm ext}$ is stable by the map $RG$.
The deviations $V^{(r,q)}[f]-V^{(r,q)}[0]$ due to the introduction of the test function $f$ now live in $\cE_{\rm ext}$.
Controlling the deviations means showing bounds on these quantities which are uniform in the UV cut-off $r$ and summable over the
scales $q$. This is crucial for the convergence of the two-sided series (\ref{dbseries}).
Let $L^{q_{+}}$ be the characteristic size of the ``support'' of $f$ (e.g., defined by the mean square distance to the origin
for a density proportional to $|f|^2$). Let $L^{-q_{-}}$ be the same notion for the Fourier transform $\widehat{f}$.
The inequality $q_{-}\le q_{+}$ is a reformulation of the uncertainty principle. The test function is approximately zero
outside a box of size $L^{q_{+}}$ and constant over cells of size $L^{q_{-}}$.
Most of the effect of the test function on the deviations and the terms in (\ref{dbseries})
concerns scales $q\in [q_{-}, q_{+}]$. One can extract (dual) exponential decay in the UV sector $q<q_{-}$ as well
as in the IR sector $q>q_{+}$. This rests on the following observations.
The fluctuation covariance $\Gamma$ decays at length scale $L$. By a suitable choice of the cut-off function
$\eta$ one can arrange for $\Gamma$ to have compact support in $x$ space (this is the idea of finite-range
decompositions in~\cite{BrydgesGM}).
If one considers the restrictions of the fluctuation field $\zeta$ to different $L$-blocks (sets of the form $L\Delta$, with
$\Delta,\Delta'$, etc. always
denoting unit cells), these are approximately independent. Thus, the RG map acts locally, i.e., can be seen
as an infinite number of independent operations performed in parallel, one for each $L$-block (the localization property).
Such an operation takes the data for the
potential from the $L^d$ unit cells $\Delta$ contained in $L\Delta'$ and produces the new data for $\Delta'$ to be used in the next RG iteration.
Another (approximate)
property of the fluctuation covariance $\Gamma$ which holds if $\eta$ is sufficiently close to $1$ near the origin is that,
almost surely, the sample
paths $\zeta$ have zero spatial average in each $L$-block. Since at the beginning of iterations $\zeta$ is smeared with the diluted test function
which is almost constant in each $L$-block, the result is zero. Thus the test function has no notable effect on the RG evolution in the UV
sector. This can also be seen in the language of Feynman diagrams where the effect of the test function is first introduced
(and later reinforced by a feedback loop) in diagrams
with external legs indicating convolutions of $\Gamma$ with the diluted test function which respectively have ``Fourier support''
in the disjoint ranges $L^{-1}\le |\xi|\le 1$ and $L^{r q_{+}}\le|\xi|\le L^{r q_{-}}$. This is a property of orthogonality
between scales (see, e.g.,~\cite{PereiraOC}).
The decay in the IR sector is due to a different mechanism: after $q_{+}-r$ RG iterations, the deviations $V^{(r,q)}[f]-V^{(r,q)}[0]$
only reside in the unit cell $\Delta_0$ at the origin. In other words, these deviations live in the span of terms
of the form $\int_{\Delta_0} :\phi^k:$ rather than $\int_{\RR^d} :\phi^k:$ which belong to the bulk.
The corresponding eigenvalues are $L^{-k[\phi]}$ instead of $L^{d-k[\phi]}$ and the RG map for the deviations becomes a contraction. 

The beauty and power of this presentation of the RG strategy is that it also works if one adds at the beginning another
term $:\phi_r^2:_{C_r}(j)$ for some new test function $j$ in order to construct the log-moment generating function $\cS^T(f,j)$
which produces mixed cumulants for both the elementary field and the squared field. 
The control of deviations in the UV sector follows a similar line of reasoning, but now is considerably more difficult.
For the BMS model, after the initial rescaling one gets
\[
:\phi_r^2:_{C_r}(j)=\int_{\RR^3} :\phi^2:_{C_0}(x)\ L^{(3-2[\phi])r}j(L^r x)\ {\rm d}^{3}x
\]
which features a new diluted form of the test function $j$.
By the localization property of the RG map and the local constancy of this diluted test function, controlling the deviations
which a priori involves the RG action in the extended space $\cE_{\rm ext}$
becomes a question which is purely about the RG action in the bulk space $\cE_{\rm bulk}$.
One has to show that if $V$ is picked as mentioned earlier on the stable manifold of $V_{\rm nontriv}$
and if $W$ is a perturbation in the $\int_{\RR^3} :\phi^2:$ direction, then (letting $n=q-r$)
the limit $\lim_{n\rightarrow\infty} RG^n(V+L^{-(3-2[\phi])n}W)$
exists and is nonzero.
However, this naive construction gives $0$ because the expanding eigenvalue at $V_{\rm nontriv}$ is strictly smaller than
the one at $V_{\rm Gauss}$, namely, $L^{3-2[\phi]}$.
Therefore the correct limit is
\[
\Psi(V,W)=\lim_{n\rightarrow\infty} RG^n(V+Z^{n} L^{-(3-2[\phi])n}W)
\]
for $Z=L^{\frac{\kappa}{2}}$ with $\kappa>0$ so that $Z^{-1}L^{3-2[\phi]}$ equals the expanding eigenvalue at $V_{\rm nontriv}$.
If $\Phi$, $\Phi^2$ respectively denote the self-similar random field and its suitably renormalized square obtained at the end of the day,
then their covariances satisfy
\[
{\rm Cov}(\Phi^2(x_1),\Phi^2(x_2))\sim \left[{\rm Cov}(\Phi(x_1),\Phi(x_2))\right]^2\times
\frac{1}{|x_1-x_2|^\kappa} \ .
\]
Thus the composite field $\cO_A=\Phi^2$ exhibits an anomalous dimension.
In fact the construction of $\Psi$ controls the deviations $V^{(r,q)}[f]-V^{(r,q)}[0]$ but is not enough for the convergence
of the series (\ref{dbseries}), in the UV sector.
Some explicit terms linear in $j$ must be extracted from the $\delta b$'s in order to secure convergence. This accounts for an
additive correction needed to define the proper $:\phi_r^2:_{C_r}(j)$ input. The correct choice is
\[
Z^{-r}\int_{\RR^3} \left(\phi_r(x)^2-L^{-2[\phi]r}\left(C_0(0)+Y\right)\right)\ j(x)\ {\rm d}^{3}x
\]
for some suitable nonuniversal constant $Y$ (whereas $Z$ is universal which here means $g$-independent).
Since $:\phi_r^2:_{C_r}(x)=\phi_r(x)^2-L^{-2[\phi]r}C_0(0)$, this is a correction to Gaussian Wick ordering
which is already needed for the definition of $\Phi^2$ when $\Phi$ is Gaussian (see~\cite{Dobrushin1,Major} for a discussion
of composite fields related to multiple stochastic integrals, in the Gaussian case).
In terms of the original $\phi^4$-type unbounded spin system $(\phi_\ux)_{\ux\in\ZZ^3}$ whose scaling limit is taken, $C_0 (0)+Y$
represents the variance of a single spin $\phi_\ux$, and $\kappa$
gives the long-distance behavior of
\begin{equation}
{\rm Cov}(\phi_{\ux}^2,\phi_{\uy}^2)\sim \frac{1}{|\ux-\uy|^{4[\phi]+\kappa}}\ .
\label{latticephi2decay}
\end{equation}
The explicit relation between the correlators and the RG dynamical system embodied in (\ref{dbseries}),
also gives a handle on questions
related to the short distance structure of these correlators: smoothness away from the big diagonal and the OPE.
The previous construction of $\cO_A=\Phi^2$ is an example of composite field renormalization.
One can also use the OPE in order to give a tautological definition of this field (in the sense of moments).
Relating the two (i.e., proving the OPE)
amounts to showing the commutation of limits $|x_1-x_2|\rightarrow 0$ and $r\rightarrow -\infty$.
The crucial ingredient $\Psi$ is reminiscent of M{\o}ller wave operators in scattering theory which intertwine free and interacting evolutions.
In the classical (rather than quantum) context of the dynamical system $RG$, such operators realize a conjugation of the nonlinear map
to its linearization at a fixed point. If $z$ denotes a curvilinear coordinate which defines the stable manifold of
$V_{\rm nontriv}$ by $z=0$ and satisfies the relation $z(RG(V))=Z^{-1}L^{3-2[\phi]} z(V)$, i.e.,
linearizes the action of $RG$ in the unstable direction,
then an easy calculation shows that $\Psi(V,W)$ is the directional derivative of $z$ at $V$ in the direction of $W$.
Coordinates such as $z$ are called nonlinear scaling fields and are the basic ingredients of Wegner's theory~\cite{Wegner} describing
the fine features of statistical mechanics systems at criticality, such as the behavior of higher composite fields.

Now is a good time to take a pause and try to answer: what is a QFT?
A safe answer is: an infinite collection of correlators $\EE[\Phi(x_1)\cdots\Phi(x_n)]$.
If one can solve the corresponding moment problem (as~\cite{CamiaGN} did for~\cite{ChelkakHI}), then one may say: a generalized random field.
As in~\cite[\S2]{Zamolodchikov}, one could also
request the secondary structure consisting of all correlators
$\EE[\cO_{A_1}(x_1)\cdots\cO_{A_n}(x_n)]$ generated from the field $\Phi$ by the OPE.
Finally, one can identify a QFT with an ideal trajectory $\cT_{\rm ideal}$. This corresponds to the modern
view in physics
which sees a QFT as a sequence of effective theories at scales $L^q$ in theory space (see~\cite{Douglas1,Douglas2}), i.e., $\cE_{\rm bulk}$.
Correlators can be recovered, via (\ref{dbseries}), as directional derivatives around the points of $\cT_{\rm ideal}$.
Yet, these directions may go out into the bigger space $\cE_{\rm ext}$.
By choosing one of the entries of $\cT_{\rm ideal}$, say the $q=0$ entry, one can parametrize such sequences or QFTs
by the unstable manifold of $V_{\rm UV}$ which typically is finite dimensional (for instance if $V_{\rm UV}=V_{\rm Gauss}$,
$[\phi]$ is canonical and $d>2$). The reparametrization obtained by choosing a different entry in the sequence (the
same as rescaling) accounts for the old pre-Wilsonian version of the RG~\cite{StueckelbergP,GellMannL}.

\section{Hierarchical models}\label{hierarchicalsec}
Since the implementation of the RG strategy is a difficult enterprise, it is useful to have simplified models
on which to test one's methods. A important example of such is that of hierarchical models.
Let $N$ be a positive integer and suppose one has a vector of centered Gaussian random variables
$(\zeta_1^{(0)},\ldots,\zeta_N^{(0)})$ whose joint law is specified by a covariance matrix $M=(M_{ij})_{1\le i,j\le N}$
whose entries sum up to zero. This implies $\sum_{i=1}^{N} \zeta_i^{(0)}=0$ almost surely.
One can make infinitely many independent copies of this vector and obtain a lattice Gaussian random field
$(\zeta_{\ux}^{(0)})_{\ux\in\LL_0}$. Here the first layer $\LL_0$ is the set $\{1,2,3,\ldots\}$
with the copies corresponding to the groups of labels $\{kN+1,\ldots,kN+N\}$, $k\ge 0$.
One can then make independent copies $(\zeta_{\ux}^{(q)})_{\ux\in\LL_q}$ of this field indexed by integers $q\ge 0$.
This introduces new layers $\LL_q$ on which one can put a geometrical structure by identifying the points of $\LL_q$
with the $N$-groups of the previous layer $\LL_{q-1}$. 
One thus obtains a singly infinite tree structure as in the following figure where $N=3$.
\[
\parbox{8cm}{
\psfrag{a}{$\LL_0$}
\psfrag{b}{$\LL_1$}
\psfrag{c}{$\LL_2$}
\raisebox{1ex}{
\includegraphics[width=8cm]{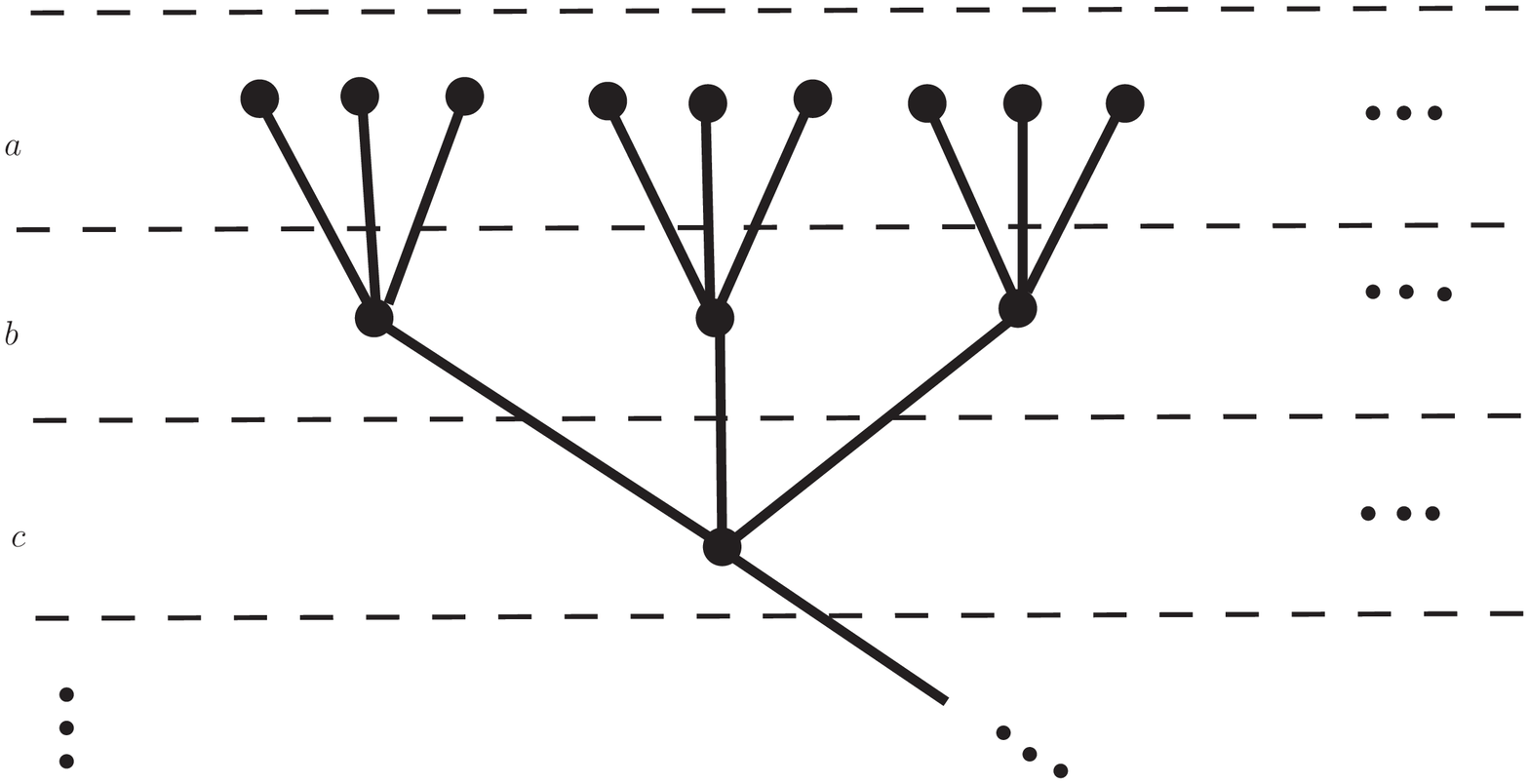}}
}
\]
For $\ux\in\LL_0$ one easily defines its ancestor $a_q(\ux)$ in $\LL_q$. Given a number $\alpha>1$,
and for $\ux\neq\uy\in\LL_0$ one defines $|\ux-\uy|=\alpha^q$ where $q$ is the smallest integer for which
$a_q(\ux)=a_q(\uy)$. This is an ultrametric notion of distance, formally denoted as the norm of a difference. 
Let $\beta>1$ be another parameter, then
\[
\phi_{\ux}=\sum_{q=0}^{\infty} \beta^{-q}\ \zeta_{a_q(\ux)}^{(q)}
\]
defines a random field $(\phi_{\ux})_{\ux\in\LL_0}$ which is a hierarchical lattice Gaussian field.
Its covariance is
$\EE[\phi_{\ux}\phi_{\uy}]\sim |\ux-\uy|^{-2[\phi]}$
with $[\phi]=\frac{\log \beta}{\log \alpha}$.
Consider $\RR^d$ discretized using dyadic cubes and identify $\LL_q$ with the set of cubes of size $2^q$.
Setting $N=2^d$, $\alpha=2$ and $\beta=2^{[\phi]}$ produces a reasonable toy model for the massless Gaussian field
on $\RR^d$ with scaling dimension $[\phi]$ and with a unit cut-off (e.g., discretized on $\ZZ^d$).
Starting from such a Gaussian hierarchical lattice measure ${\rm d}\mu_{C_0}$ for $(\phi_{\ux})_{\ux\in\LL_0}$
one can perturb it by a product of single spin potentials involving $\phi_{\ux}^{2}$ and $\phi_{\ux}^{4}$ terms and repeat the previous story
by integrating out the fluctuation fields $\zeta$ a few layers at a time.
Namely, one can use the RG strategy in order to analyze the scaling limit of the resulting non-Gaussian field.
A nice feature of such toy models is that many of the properties which earlier were approximately true, for instance
the localization property of the RG, now become exact.
In order to have a home for the scaling limit, one needs a notion of continuum. This is obtained by subdividing the nodes of the top layer
$\LL_0$ and 
continuing the tree structure with the introduction of new layers $\LL_{-1}$, $\LL_{-2}$, etc.
The set $\LL_{-\infty}$ of leafs or ends at infinity of the resulting doubly infinite tree structure is the needed continuum.
The lattice $\LL_0$ is to $\LL_{-\infty}$ what the lattice $\ZZ^d$ is to $\RR^d$.

The introduction of hierarchical models originated from two independent sources.
One is the work of Dyson on one-dimensional Ising spin models
with long-range interactions~\cite{Dyson}. The other is the early work of Wilson on his RG theory.
While some features of the latter were already present in the old article~\cite{Wilson65}, one may say that its first systematic exposition
was given in~\cite{WilsonII}. In fact, this article was about a hierarchical model as above and the corresponding RG map was called
``the approximate recursion''. The relevance of this approximation for physical models over $\RR^d$ stemmed from
Wilson's anticipation, in that article, of wavelet multiresolution analysis (see~\cite{Battle}).
Note that the word ``approximation'' for the hierarchical RG in relation to the RG for real models is somewhat misleading, since
it suggests that this is the zero-th step of a systematic approximation procedure which can lead, through successive improvements,
to the ``true'' RG for models over $\RR^d$. Despite some attempts in this direction (e.g.,~\cite{MeuriceApx}
and~\cite[\S14.2]{Meurice}), the existence of such a procedure is unclear.
The confusion may be due to the fact the hierarchical RG is thematically similar to another simplification called the local potential
approximation (LPA). The LPA arises from the above description of the RG on $\RR^d$ if one ignores the nonlocal kernels
which should appear in $V'=RG(V)$ as well as the gradient terms one obtains by the comparison of such nonlocal kernels with their
local projection. 
There is rigorous work in the LPA setting, e.g.,~\cite{Felder} as well as a nonrigorous approximation procedure
starting from the LPA know as the derivative expansion (see, e.g.,~\cite{BagnulsB}).
The author's point of view is that one should not try to approximate critical exponents such as anomalous dimensions
for models on $\RR^d$ by their analogues on hierarchical models. The latter reflect the different tree-like geometric texture
of the underlying continuum (and in particular depend on the choice of $N$ which governs the shape of the underlying tree-like space).
The utility of hierarchical models is that they are a good
testing ground for RG techniques. Such methodology has been successful, in a rigorous setting, for instance in the work of
Gaw\c{e}dzki and Kupiainen on $(\nabla\phi)^4$ lattice models (the hierarchical model testing was done in~\cite{GawedzkiKgr1,GawedzkiKgr2}
while the real model was treated in~\cite{GawedzkiKgr3,GawedzkiKgr4}) or that of Brydges
and Slade~\cite{BrydgesS} on the weakly self-avoiding walk in
four dimensions (their approach was tested on a hierarchical model in~\cite{BrydgesI1,BrydgesI2}).
Another example, in a nonrigorous context, of the success of this methodology is the work of Wilson himself
when he developed his RG theory in the first place.
Indeed, this theory presented in~\cite{WilsonII} initially had a modest impact, perhaps because it pertained
to the hierarchical model and it was not clear, to the skeptical minds of the time, how it could shed light
on the physics of real models. The situation drastically changed with the soon-to-follow article~\cite{WilsonF}
which allowed the RG method to produce an expansion (which can be systematically improved) for critical
exponents of real models such as the Ising model in three dimensions related to the liquid-vapour critical point of water
(i.e., a real-world phenomenon where such exponents can be measured experimentally).
This so-called $\epsilon$-expansion whose definitive treatment was later given in~\cite{WilsonK} was largely responsible for
the revolution created by Wilson's RG theory in physics. It is now part of any theoretical
physicist's DNA or view of the world (see, e.g.,~\cite{Douglas1,Douglas2}).
The key article~\cite{WilsonF} contained two major conceptual advances. The first is that of introducing the bifurcation
parameter $\epsilon$ (as done later in the BMS model) and expanding exponents with respect to this parameter which can be viewed as
the difference $4-d$ between spatial dimensions. This idea can be implemented for both the $\RR^d$ model
and the hierarchical one. The second idea was the implementation of this $\epsilon$-deformation in the $\RR^d$ case using
the newly developed dimensional regularization in QFT. One can thus ask if the first idea was initially developed on the hierarchical model.
The answer is ``yes'' or in Wilson's words ``Then, at Michael's urging,
I work out what happens near four dimensions for the approximate
recursion formula, and find that d-4 acts as a small parameter. Knowing this it is then trivial,
given my field theoretic training, to construct the beginning of the epsilon expansion for
critical exponents.''~\cite{WilsonInterview}.

There is great arbitrariness when setting up a hierarchical model in order to mimic one living in $\RR^d$.
More precisely, there are lots of ways to pick $N$, $M$, $\alpha$ and $\beta$ for given $d$ and $[\phi]$.
There are many versions of the hierarchical model considered by various authors (see~\cite{Meurice} for a review).
If one does not set up this model carefully, one can end up with rather absurd results such as lack of universality
or having critical exponents produced by the RG strategy depend on the artificial yardstick $L$
(see, e.g.,~\cite[\S5.2]{Meurice} for a discussion of this issue).
There is a particular set-up for the hierarchical RG which, among many other beautiful mathematical properties, avoids such problems
and it uses $p$-adics. For $d$ an integer,
the $p$-adic set-up consists in taking $N=p^d$ where $p$ is a prime number, $\alpha=p$, $\beta=p^{[\phi]}$
and defining the matrix $M$ by putting $1-p^{-d}$ on the diagonal and $-p^{-d}$ everywhere else.
The $p$-adic BMS model is the particular case of hierarchical model obtained in this way when $d=3$ and
$[\phi]=\frac{3-\epsilon}{4}$.
From the point of view of probability theory, the restriction to primes is not essential. However, doing so gives access
to a huge ``software library'' developed for the needs of number theory.
One can then identify $\LL_{-\infty}$ with $\QQ_p^d$ where $\QQ_p$ is the field (in the algebra sense) of $p$-adic numbers.
The fields (in the QFT or probability theory sense) are still real-valued.
Instead of the previous elementary and ad hoc description of the $p$-adic model, a more elegant approach (see~\cite{AbdesselamCG1})
is to use Fourier analysis and the theory of distributions on the $p$-adics, for real or complex-valued functions,
which come from this software library. Since $\QQ_p$ is an additive group, there is a natural notion of translation invariance.
The maximal compact subgroup $GL_d(\ZZ_p)$ (unique up to conjugation), which is the analogue of the orthogonal group in $\RR^d$,
supplies the notion of rotation invariance. The analogue of the Euclidean norm in $\RR^d$ is the maximum of the
$p$-adic absolute values
of the components (used to define the previous hierarchical distance $|\ux-\uy|$) since it is invariant by $GL_d(\ZZ_p)$.
Instead of $\RR_{+}^{\ast}$, the group of scaling transformations $p^{\ZZ}$
is now discrete.
One sets $L=p^l$ for some integer $l\ge 1$ when defining the RG which then integrates out the first $l$ layers $\LL_0,\ldots,\LL_{l-1}$
at each step.
One thus avoids the problem of
$L$-dependent critical exponents (since the texture of the underlying space depends on $p$, not $L$).
This is also the same RG map discussed earlier using a Fourier cut-off.
Indeed, by taking the function $\eta$ equal to the sharp characteristic function of the interval $[0,1]$
one recovers the above ad hoc hierarchical model.
One can also try to localize the previous notions of invariances and look for an analogue of conformal
invariance (see~\cite{Melzer,LernerM,Lerner}).
A leitmotiv in the theory of automorphic forms is that completions of $\QQ$
such as $\RR$ or $\QQ_p$ should be treated on equal footing. Beautiful theories in analysis
(e.g., the theory of unitary representations of noncompact
groups) developed for $\RR$ have analogues over $\QQ_p$. Such unity is visible in the series of books by Gel'fand and co-authors
on generalized functions (e.g.,~\cite{Gelfand1,Gelfand4} which pertain to self-similar random fields)
which included~\cite{Gelfand6} as a sixth volume in the Russian edition. 
The author believes (and hopes to have convinced the reader) that a similar harmonious unity is present
in the context of generalized random fields, QFT and the RG.

The previous presentation of QFT and the RG may seem somewhat impressionistic and it indeed avoided discussing many issues
which are important for mathematical rigor: the infinite volume limit, dealing with nonlocalities and gradient terms,
the specific norms and bounds needed, etc.
Nevertheless, the reader can find in the article~\cite{AbdesselamCG1} a complete rigorous substantiation of the story
told in \S\ref{fundamentalsec} and \S\ref{strategysec} (with the exception of the OPE),
in the case of the self-similar $p$-adic BMS model at $V_{\rm nontriv}$.
For the model over $\RR^3$, the reader is referred to the preliminary results~\cite{BMS,AbdesselamCMP} which should be easier to read after
seeing a simpler version~\cite[\S6]{AbdesselamCG1} of the needed RG estimates.

\bigskip
\noindent{\bf Acknowledgements:}
{\small
The author would like to express his gratitude to those who influenced his thinking about QFT and the RG, over the years.
These are D. C. Brydges, J. Magnen, P. K. Mitter and V. Rivasseau. Of course, any shortcoming of the present article is the responsibility
of the author alone.}


\begin{thebibliography}{999}

\bibitem{AbdesselamCMP}
A.~Abdesselam,
A complete renormalization group trajectory between two fixed points.
Comm. Math. Phys. {\bf 276} (2007), no. 3, 727–-772. 

\bibitem{AbdesselamCG1}
A. Abdesselam, A. Chandra and G. Guadagni,
Rigorous quantum field theory functional integrals over the $p$-adics I: anomalous dimensions.
Preprint arXiv:1302.5971[math.PR], 2013.

\bibitem{ATLAS}
ATLAS Collaboration,
Observation of a new particle in the search for the Standard Model Higgs boson with the ATLAS detector at the LHC.
Phys. Lett. B {\bf 716} (2012), no. 1, 1--29.

\bibitem{BagnulsB}
C. Bagnuls and C. Bervillier,
Exact renormalization group equations: an introductory review.
Renormalization group theory in the new millennium, II.
Phys. Rep. {\bf 348} (2001), no. 1--2, 91–-157. 

\bibitem{Battle}
G. Battle,
Wavelet refinement of the Wilson recursion formula. In: ``{\em Recent Advances in Wavelet Analysis}'', 87-–118,
Wavelet Anal. Appl., {\bf 3}, Academic Press, Boston, MA, 1994. 


\bibitem{BeilinsonD}
A.~Beilinson and V.~Drinfeld,
{\em Chiral Algebras}.
American Mathematical Society Colloquium Publications, 51. American Math. Soc., Providence, RI, 2004.


\bibitem{Borcherds}
R.~E.~Borcherds,
Vertex algebras, Kac-Moody algebras, and the monster.
Proc. Nat. Acad. Sci. U.S.A. {\bf 83} (1986), no. 10, 3068-–3071. 

\bibitem{BrydgesGM}
D. C. Brydges, G. Guadagni, and P. K. Mitter,
Finite range decomposition of Gaussian processes.
J. Statist. Phys. 115 (2004), no. 1-2, 415–449. 

\bibitem{BrydgesI1}
D. C. Brydges and J. Z. Imbrie,
End-to-end distance from the Green's function for a hierarchical self-avoiding walk in four dimensions.
Comm. Math. Phys. {\bf 239} (2003), no. 3, 523–-547. 

\bibitem{BrydgesI2}
D. C. Brydges and J. Z. Imbrie,
Green's function for a hierarchical self-avoiding walk in four dimensions.
Comm. Math. Phys. {\bf 239} (2003), no. 3, 549-–584. 

\bibitem{BMS}
D.~C.~Brydges, P.~K.~Mitter and B.~Scoppola.
Critical $(\Phi^4)_{3,\epsilon}$.
Comm.~Math.~Phys., {\bf 240} (2003), 281--327.

\bibitem{BrydgesS}
D. Brydges and G. Slade,
Renormalisation group analysis of weakly self-avoiding walk in dimensions four and higher. In: ``{\em Proceedings of the International
Congress of Mathematicians}'', Vol. IV, 2232–-2257, Hindustan Book Agency, New Delhi, 2010. 

\bibitem{CamiaGN}
F. Camia, C. Garban and C. Newman,
Planar Ising magnetization field I. Uniqueness of the critical scaling limit.
Preprint arXiv:1205.6610[math.PR], 2012, to appear in Ann. Probab.

\bibitem{CandelasOGP}
P. Candelas, X. C. de la Ossa, P. S. Green and L. Parkes,
A pair of Calabi-Yau manifolds as an exactly soluble superconformal theory.
Nuclear Phys. B {\bf 359} (1991), no. 1, 21–-74. 

\bibitem{ChelkakHI}
D. Chelkak, C. Hongler and K. Izyurov,
Conformal invariance of spin correlations in the planar Ising model.
Preprint arXiv:1202.2838[math-ph], 2012.

\bibitem{CMS}
CMS Collaboration,
Observation of a new boson at a mass of 125 GeV with the CMS experiment at the LHC.
Phys. Lett. B {\bf 716} (2012), no. 1, 30--61.


\bibitem{Costello}
K. Costello,
{\em Renormalization and effective field theory}.
Mathematical Surveys and Monographs, {\bf 170}. American Mathematical Society, Providence, RI, 2011. 

\bibitem{CostelloG}
K. Costello and O. Gwilliam,
{\em Factorization Algebras in Quantum Field Theory}, book in progress.
Available at \verb+http://www.math.northwestern.edu/~costello/factorization.pdf+

\bibitem{Dobrushin1}
R. L. Dobrushin,
Gaussian and their subordinated self-similar random generalized fields.
Ann. Probab. {\bf 7} (1979), no. 1, 1-–28. 

\bibitem{Dobrushin2}
R. L. Dobrushin,
Automodel generalized random fields and their renorm group. In: {\it Multicomponent Random Systems}, Ed.: R. L. Dobrushin and
Ya. G. Sinai, pp. 153-–198,
Adv. Probab. Related Topics {\bf 6}, Marcel Dekker, New York, 1980. 

\bibitem{Douglas1}
M. R. Douglas,
Spaces of quantum field theories.
Preprint arXiv:1005.2779[hep-th], 2010.

\bibitem{Douglas2}
M. R. Douglas,
Foundations of quantum field theory. In: ``{\em String-Math 2011}'', 105-–124,
Proc. Sympos. Pure Math., {\bf 85}, American Math. Soc., Providence, RI, 2012. 

\bibitem{Dyson}
F. J. Dyson,
Existence of a phase-transition in a one-dimensional Ising ferromagnet.
Comm. Math. Phys. {\bf 12} (1969), no. 2, 91–-107. 

\bibitem{deFariaM}
E. de Faria and W. de Melo,
{\em Mathematical aspects of quantum field theory}.
Cambridge Studies in Advanced Mathematics, {\bf 127}. Cambridge University Press, Cambridge, 2010. 

\bibitem{Felder}
G. Felder,
Renormalization group in the local potential approximation.
Comm. Math. Phys. {\bf 111} (1987), no. 1, 101–-121. 

\bibitem{Folland}
G. B. Folland,
{\em Quantum field theory.
A tourist guide for mathematicians}. Mathematical Surveys and Monographs, {\bf 149}. American Mathematical Society, Providence, RI, 2008.

\bibitem{Frenkel}
E. Frenkel,
{\em Langlands correspondence for loop groups}.
Cambridge Studies in Advanced Mathematics, {\bf 103}. Cambridge University Press, Cambridge, 2007. 

\bibitem{FrenkelB}
E. Frenkel and D. Ben-Zvi,
{\em Vertex Algebras and Algebraic Curves}.
Mathematical Surveys and Monographs, {\bf 88}. American Math. Soc., Providence, RI, 2001. 

\bibitem{FultonRW}
T. Fulton, F. Rohrlich and L. Witten,
Conformal Invariance in Physics.
Rev. Mod. Phys. {\bf 34} (1962), no. 3, 442-–457.

\bibitem{Gaitsgory}
D.~Gaitsgory,
Notes on 2D conformal field theory and string theory. In {\em Quantum Fields and Strings: a Course for Mathematicians},
Vol. 2 (Princeton, NJ, 1996/1997), Edited by P. Deligne et al.
pp. 1017-–1089, American Math. Soc., Providence, RI, 1999. 

\bibitem{GawedzkiKgr1}
K. Gaw\c{e}dzki and A. Kupiainen,
Renormalization group study of a critical lattice model. I. Convergence to the line of fixed points.
Comm. Math. Phys. {\bf 82} (1981/82), no. 3, 407–-433. 

\bibitem{GawedzkiKgr2}
K. Gaw\c{e}dzki and A. Kupiainen,
Renormalization group study of a critical lattice model. II. The correlation functions.
Comm. Math. Phys. {\bf 83} (1982), no. 4, 469-–492. 

\bibitem{GawedzkiKgr3}
K. Gaw\c{e}dzki and A. Kupiainen,
Block spin renormalization group for dipole gas and $(\nabla\varphi)^4$.
Ann. Physics {\bf 147} (1983), no. 1, 198-–243. 

\bibitem{GawedzkiKgr4}
K. Gaw\c{e}dzki and A. Kupiainen,
Lattice dipole gas and $(\nabla\varphi)^4$ models at long distances: decay of correlations and scaling limit.
Comm. Math. Phys. {\bf 92} (1984), no. 4, 531–-553. 


\bibitem{Gelfand6}
I. M. Gel'fand, M. I. Graev, M. I. and I. I. Pyatetskii-Shapiro,
{\em Representation Theory and Automorphic Functions}.
Translated by K. A. Hirsch. W. B. Saunders Co., Philadelphia--London--Toronto, 1969.

\bibitem{Gelfand1}
I. M. Gel'fand and G. E. Shilov,
{\em Generalized Functions. Vol. 1: Properties and Operations}.
Translated by E. Saletan. Academic Press, New York--London, 1964.

\bibitem{Gelfand4}
I. M. Gel'fand and N. Ya. Vilenkin,
{\em Generalized Functions. Vol. 4: Applications of Harmonic Analysis}.
Translated by A. Feinstein. Academic Press, New York--London, 1964. 

\bibitem{GellMannL}
M. Gell-Mann and F. E. Low,
Quantum electrodynamics at small distances. 
Phys. Rev. (2) {\bf 95} (1954), 1300--1312.


\bibitem{KangM}
N.-G.~Kang, and N.~Makarov,
{\em Gaussian Free Field and Conformal Field Theory}.
Ast\'erisque {\bf 353}. Soc. Math. de France, 2013.
Also available as preprint arXiv:1101.1024v2[math.PR], 2011.

\bibitem{KapustinW}
A. Kapustin and E. Witten,
Electric-magnetic duality and the geometric Langlands program.
Commun. Number Theory Phys. {\bf 1} (2007), no. 1, 1-–236. 

\bibitem{Lerner}
\`E. Yu. Lerner,
The hierarchical Dyson model and $p$-adic conformal invariance.
Theor. Math. Phys. {\bf 97} (1993), no. 2, 1259–-1266.

\bibitem{LernerM}
\`E. Yu. Lerner and M. D. Missarov,
$p$-adic conformal invariance and the Bruhat-Tits tree.
Lett. Math. Phys. {\bf 22} (1991), no. 2, 123-–129. 

\bibitem{Lurie}
J. Lurie,
On the classification of topological field theories. Current developments in mathematics, 2008, 129-–280, Int. Press, Somerville, MA, 2009. 

\bibitem{Major}
P. Major,
{\it Multiple Wiener-It\^{o} Integrals. With Applications to Limit Theorems}.
Lecture Notes in Mathematics {\bf 849}, Springer, Berlin, 1981. 

\bibitem{Melzer}
E. Melzer,
Non-Archimedean conformal field theories.
Internat. J. Modern Phys. A {\bf 4} (1989), no. 18, 4877–-4908. 

\bibitem{MeuriceApx}
Y. Meurice,
A perturbative improvement of the hierarchical approximation.
Unpublished preprint arXiv:hep-th/9307128, 1993.

\bibitem{Meurice}
Y. Meurice,
Nonlinear aspects of the renormalization group flows of Dyson's hierarchical model.
J. Phys. A {\bf 40} (2007), no. 23, R39-–R102. 


\bibitem{PereiraOC}
E. Pereira and M. O'Carroll,
Orthogonality between scales and wavelets in a representation for correlation functions.
The lattice dipole gas and $(\nabla\phi)^4$ models.
J. Statist. Phys. {\bf 73} (1993), no. 3--4, 695–-721. 

\bibitem{Polyakov}
A. M. Polyakov,
Conformal symmetry of critical fluctuations.
J. Exp. Theor. Phys. Lett. {\bf 12} (1970), 381--383.

\bibitem{Schramm}
O. Schramm,
Conformally invariant scaling limits: an overview and a collection of problems. In: {\em International Congress of
Mathematicians}, Vol. I, 513–-543, European Math. Soc., Z\"{u}rich, 2007. 

\bibitem{Smirnov}
S. Smirnov,
Discrete complex analysis and probability. In: ``{\em Proceedings of the International Congress of Mathematicians}'', Vol. I,
595-–621, Hindustan Book Agency, New Delhi, 2010. 

\bibitem{StueckelbergP}
E. C. G. Stueckelberg and A. Petermann,
La normalisation des constantes dans la th\'eorie des quanta.
Helvetica Phys. Acta {\bf 26} (1953), 499--520.


\bibitem{Wegner}
F. J. Wegner,
Corrections to scaling laws.
Phys. Rev. B {\bf 5} (1972), no. 11, 4529-–4536.

\bibitem{Wilson65}
K. G. Wilson,
Model Hamiltonians for local quantum field theory.
Phys. Rev. {\bf 140} (1965), no. 2B, B445-–B457. 

\bibitem{WilsonII}
K. G. Wilson,
Renormalization group and critical phenomena. II. Phase-space cell analysis of critical behavior.
Phys. Rev. B {\bf 4} (1971), no. 9, 3184-–3205. 

\bibitem{WilsonF}
K. G. Wilson and M. E. Fisher,
Critical Exponents in 3.99 Dimensions.
Phys. Rev. Lett {\bf 28} (1972), no. 4, 240--243.

\bibitem{WilsonK}
K. G. Wilson and J. Kogut,
The renormalization group and the $\epsilon$ expansion.
Phys. Rep.  {\bf 12} (1974), no. 2, 75–-199.

\bibitem{WilsonInterview}
K. G. Wilson, cited from Part II of his 07/06/2002 interview in Physics of Scales Activities.
Transcript available at
\verb+http://authors.library.caltech.edu/5456/1/hrst.mit.edu/hrs/renormalization/Wilson/Wilson2.htm+

\bibitem{WittenKnots}
E. Witten,
Quantum field theory and the Jones polynomial.
Comm. Math. Phys. {\bf 121} (1989), no. 3, 351-–399. 

\bibitem{Witten4manifolds}
E. Witten,
Monopoles and four-manifolds.
Math. Res. Lett. {\bf 1} (1994), no. 6, 769–-796. 

\bibitem{WittenOPE}
E. Witten,
Perturbative quantum field theory. 
In {\em Quantum Fields and Strings: a Course for Mathematicians},
Vol. 1 (Princeton, NJ, 1996/1997), Edited by P. Deligne et al.
pp.  419–-473, American Math. Soc., Providence, RI, 1999. 

\bibitem{Wu}
T. T. Wu,
Theory of Toeplitz determinants and the spin correlations of the two-dimensional Ising model. I.
Phys. Rev. {\bf 149} (1966), no. 1, 380-–401.


\bibitem{Zamolodchikov}
A.~B.~Zamolodchikov,
Renormalization group and perturbation theory about fixed points in two-dimensional field theory.
Sov.~J.~Nucl.~Phys. {\bf 46} (1987), 1090–-1096. 

\end{thebibliography}
\end{document}